\documentclass[11pt]{article}

\usepackage{latexsym}
\usepackage{amssymb}
\usepackage{amsmath}
\usepackage{epsfig}
\usepackage{multicol}
\usepackage{graphicx}

\newtheorem{theorem}{Theorem}[section]

\newtheorem{proposition}[theorem]{\bf{Proposition}}

\newtheorem{remark}[theorem]{\bf{Remark}}

\begin{document}
%
\title{On a Topological Problem of Strange Attractors}
%
%
%
\author{
 Ibrahim Kirat and Ayhan Yurdaer \\
\footnotesize Department of Mathematics, Istanbul Technical
University, 34469,Maslak-Istanbul, Turkey \\
\footnotesize E-mail: ibkst@yahoo.com \ and \ yurdaerayhan@itu.edu.tr }

%

%
%

\date {  }
\maketitle

\begin{abstract}
Somehow, the revised version of our paper \cite{KY} does not appear on journals' home page. Here we present the revised version altered to reflect the corrections
and/or additions to that paper. In this note, we consider self-affine attractors that are generated by an integer expanding $n\times n$ matrix (i.e., all of its eigenvalues have moduli $>1$) and a finite set of vectors in  ${\Bbb{Z}}^n$. We concentrate on the problem of connectedness for $n\leq 2$. Although, there has been intensive study on the topic recently, this problem is not settled even in the one-dimensional case. We focus on some basic attractors, which have not been studied fully, and characterize connectedness.

\bigskip

\end{abstract}

\section{Introduction}\label{INTRO}

Let $S_1,...,S_q$, $q>1$, be contractions on ${\Bbb{R}}^n$, i.e.,
$||S_j(x)-S_j(y)||\leq c_j||x-y||$ for all $x, y \in {\Bbb{R}}^n$ with
$0<c_j<1$. Here $||\cdot ||$ stands for the usual Euclidean norm, but this norm may be replaced by any other norm on ${\Bbb{R}}^n$. It is well known \cite{F2} that there exists a unique
non-empty compact set $F\subset {\Bbb{R}}^n$ such that
$$
 F = \bigcup_{j=1}^q S_j(F).
$$

Let $M_n({\Bbb R})$ denote the set of $n\times n$ matrices
with real entries. We will assume that
$$
S_j(x)=T^{-1}(x+d_j),  \ \ \ \ x\in {\Bbb{R}}^n,
$$
where $d_j \in {\Bbb{R}}^n$, called \emph{digits}, and $T \in  M_n({\Bbb R}).$
Then $F$ is called a \emph{self-affine set} or a \emph{self-affine fractal},
and can be viewed as the invariant set or the attractor of the (affine)
{\it iterated function system} (IFS)
$\{S_j(x)\}$ (in the terminology of dynamical systems). Let $M_n({\Bbb{Z}})$ be the set of $n\times n$ integer matrices.
Further, if
$D:=\{d_1,...,d_q\}\subset {\Bbb{Z}^n}$ and
$T\in M_n({\Bbb{Z}})$, it is called an \emph{integral self-affine set}
and we will primarily consider such sets in this paper.
 If, additionally, $|\det(T)|=q$
and the integral self-affine set $F$ has positive Lebesgue measure, then $F$
is called an \emph{integral self-affine tile}. We sometimes write $F(T, D)$ for $F$ to
stress the dependence on $T$ and $D$. For such tiles, the positivity of the Lebesgue measure
is equivalent to having nonempty interior \cite{B}.

There is a demand to develop
analysis on fractal spaces, in order to deal with physical phenomena like heat and electricity flow in disordered media, vibrations of fractal materials and turbulence in fluids.
Without a better understanding of the topology of fractals, this seems to be a difficult task. There is a growing literature on the formalization and
representation of topological questions; see \cite{DEG} for a survey
of the field.

One of the interesting aspects of the self-affine sets
is the \textit{connectedness}, which roughly means the attractor cannot be written as a disjoint union of two pieces. This property is important in computer vision and remote sensing \cite{HZ,R}. We mention that connected self-affine fractals
are curves; thus, they are sometimes referred to as self-affine curves \cite{K2}. There is
some motivation for studying connected self-affine tiles because
they are related to number systems, wavelets, torus maps.  Recently,
there have been intensive investigations on the topic by Kirat and Lau \cite{KL1,K2}, Akiyama and
Thuswaldner \cite{AT,LAT}, Ngai and Tang \cite{NT1,NT2} and Luo et
al. \cite{LAT,LRT}.


In this note, we consider planar integral self-affine fractals obtained from $2\times 2$ integer
matrices with reducible characteristic polynomials, and report our findings on their connectedness.
However, our considerations can be generalized to higher dimensions. As for the organization of the
paper, in Section \ref{simple}, we deal with special cases
and state some simple, but unconventional techniques to check the connectedness.
In Section \ref{GENERAL}, we study the neighbor sets of self-affine fractals.

\section{Some Unconventional Techniques }\label{simple}

Usually, connectedness criteria were given by using a
``graph'' with vertices in $D$ \cite{H,KL1}. In this section,
we present graph-independent techniques to check the
connectedness or disconnectedness.
Throughout the paper, $T^{-1}$ is a contraction. Let $\# D$ denote the number of elements in $D$. We first recall a known result.

\begin{proposition}\label{collinear1} \cite{KL1} \
Suppose $T=[\pm q] $ with $q \in {\Bbb N}$, and
$D \subseteq {\Bbb R}$ with $\# D=q.$ Then $F(T,D)$ is a connected tile
if and only if, up to a translation,
$D=\{0,a,2a,...,(q-1)a \}$ for some $a>0$.
\end{proposition}

As one may notice $q$ and $D$ are not arbitrary in Proposition \ref{collinear1} since $q \in {\Bbb N}$ and $\# D=q$.
By using the approach in \cite{K1,KK1}, we can remove such restrictions. For that purpose, we consider the convex hull of
$F$ and denote it by $K$. Also let $K_1=\bigcup_{j=1}^q S_j(K).$ Then we have the following.

\begin{proposition}\label{collinear2} \ Let $D=\{0,d_2v,\cdots, d_{q} v\}\subset {\Bbb{R}^n}$ with $v\in \mathbb{R}^n\setminus \{0\}$
and $T= p I$, where  $p\in \mathbb{R}$ and $I$ is the identity matrix. Then $F(T,D)$ is connected if and only if $K=K_1$.
\end{proposition}

\begin{figure}[!hbt]\label{fig1}
\centering
\scalebox{0.5}{\includegraphics*{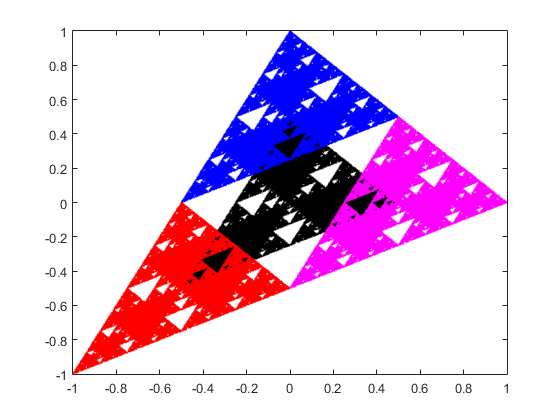}}
\caption{The Sierpi$\acute{\rm n}$ski tile}\label{fig1}
\end{figure}



\begin{remark} {\rm
 A digit set $D$ as in Proposition \ref{collinear2} is called a \textit{collinear} digit set. It is easy to check the condition $K=K_1$ in the proposition because $K$ is a closed interval. Also, note that if $T=\pm 2 I$, then $F(T,D)$ is connected for any digit set. A famous example of this
type is the Sierpi$\acute{\rm n}$ski tile (see Figure \ref{fig1}),  for which $T= 2 I$ and
$D=\{ \tiny d_1=\left[
\begin{array}{cc}  0 \\ 0  \end{array}
\right], d_2=\left[
\begin{array}{cc}  1 \\ 0  \end{array}
\right], d_3=\left[
\begin{array}{cc}  0 \\ 1  \end{array}
\right],
 d_{4}=\left[
\begin{array}{cc}  -1 \\ -1  \end{array}
\right]
 \}$. }

\end{remark}

The disconnectedness of $F(T,D)$ was studied in \cite{Ko}.
Here we want to mention another unconventional sufficient condition for disconnectedness.
In the rest of the paper, we study attractors $F(T,D)$ in the plane  such that $T\in M_2({\Bbb Z})$ has a reducible characteristic polynomial. From \cite{K2}, we know that such matrices are conjugate to one of the following lower triangular matrices
\begin{equation}\label{matrix}
{\footnotesize \left[
\begin{array}{cc} n  & 0 \\ t & m \end{array} \right]}, \ \ \ \  \mbox{where $|n| \geq |m|$, and $t=0$ or $t=1$} .
\end{equation}
We also let
$$S= \{ {\tiny \left[
\begin{array}{cc} i  \\ j  \end{array} \right]} :  0\leq i \leq |n| -1, \ 0\leq j \leq |m| -1\}.$$
 The attractors of the next proposition can be considered as a generalization of Sierpi$\acute{\rm n}$ski carpets \cite{M}.
Let $dim_S(F)$ be the \emph{singular value dimension} of $F$ (see \cite{F3}). We call a collinear digit set $D$ with $v$ is an eigenvector of $T$ \textit{eigen-collinear}.
In that case, $F$ is a subset of a line segment. By using Corollary 5 in \cite{F3}, we obtain the following.

\begin{proposition}\label{disconnectedness} \ Assume that $T$ is as in (\ref{matrix}), $D\subset S$, and $D$
is not eigen-collinear. Then $F(T,D)$ is disconnected
if $\log_{|m|} r+\log_{|n|} (\frac{q}{r})\neq dim_S(F)$, where $q=\# D$ and $r$ is the number of $j$ so that ${\tiny \left[
\begin{array}{cc} i  \\ j  \end{array} \right]}\in D$ for some $i$.
\end{proposition}

\begin{remark} {\rm It is easy to check the sufficient condition for the attractors $F(T,D)$ in Proposition \ref{disconnectedness} because, in that case,
\[ dim_S(F)=\left\{
\begin{array}{cc}
             1+\log_{|n|}(\frac{q}{|m|})  & \ \ \  \mbox{if $|m| <q\leq |mn|$,} \\
              \log_{|m|} q & \mbox{if  $q \leq |m|$.}
             \end{array} \right. \] }

\end{remark}

\section{On the Neighbor Sets}\label{GENERAL}

In this section, we will present a practical way of checking the connectedness of $F(T,D)$ with $T$ as in (\ref{matrix}) and $D\subset S$. Note that it is enough to consider the case $n,m>0$, since $F(T,D)=F(T^2,D+TD)$.
Let $\mathcal{N}=(F-F)\cap ({\Bbb Z}^2\setminus \{\tiny \left[\begin{array}{cc} 0  \\ 0  \end{array} \right]\})$, which we call the \emph{neighbor set} of $F$.
Set
$$
 \Delta D=D-D, \ \ \
 {
a_1=\tiny \left[\begin{array}{cc} n-1  \\ 1  \end{array} \right], \ \ \
a_2=\tiny \left[\begin{array}{cc} 0  \\ m-1  \end{array} \right], \ \ \
a_3=\tiny \left[\begin{array}{cc} n-1  \\ m-1  \end{array} \right]}, \ \ \
b_1=\tiny \left[\begin{array}{cc} n-1  \\ 0  \end{array} \right],
$$
$$
{e_1=\tiny \left[\begin{array}{cc} 1  \\ 0  \end{array} \right], \ \
e_2=\tiny \left[\begin{array}{cc} 0  \\ 1  \end{array} \right], \ \
e_3=\tiny \left[\begin{array}{cc} 1  \\ 1  \end{array} \right], \ \
e_4=\tiny \left[\begin{array}{cc} 1  \\ -1  \end{array} \right] }.
$$

\begin{proposition}\label{neighbors1} \ Assume that $F$ is as in Proposition \ref{disconnectedness}, $t=1$ and $n,m>0$. Then

(i)  if $ a_2 \notin \Delta D$, then $F$ is disconnected,

(ii) otherwise,

 $\mathcal{N}= \{ \pm e_i \ | \ i\in \{1,2 \} \ {\rm and} \ a_i \in \Delta D \}\cup \{ \pm e_4 \ | \  \ a_1- a_2 \in \Delta D \} $.
\end{proposition}

\begin{proposition}\label{neighbors2} \ Assume that $F$ is as in Proposition \ref{disconnectedness}, $t=0$ and $n,m>0$. Let $b_2=a_2$, $b_3=a_3$.  Then

(i)  if $ b_1, b_2, b_3, b_1- b_2  \notin \Delta D$, then $F$ is disconnected,

(ii) otherwise,

 $\mathcal{N}= \{ \pm e_i \ | \ i\in \{1,2,3 \} \ {\rm and} \ b_i \in \Delta D  \}\cup \{ \pm e_4 \ | \  \ b_1-  b_2 \in \Delta D \}$.
\end{proposition}

For a digit set $D$,
 an {\it s-chain} (in $D$) is a finite sequence $\{d_1,...,d_s \}$ of $s$
 vectors in $D$ such that $d_i-d_{i+1}\in \mathcal{N}$ for $i=1,...,s-1.$ Then we can put the connectedness criterion in \cite{KL1} into the following form.


\begin{proposition}\label{connectedness} \  $F$ is connected if and only if, by re-indexing $D$ (if necessary),
$D$ forms a $q$-chain.
\end{proposition}

\begin{figure}[!hbt]\label{fig2}
\centering
\scalebox{0.3}{\includegraphics*{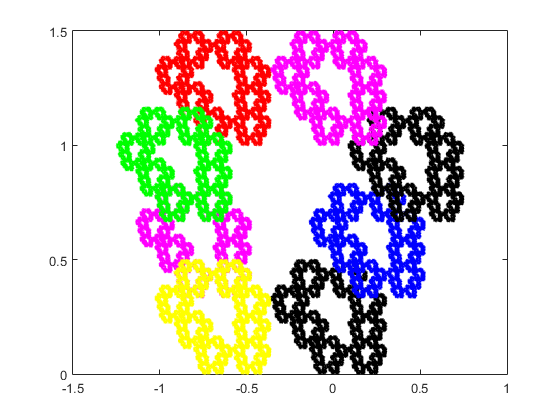}}\hspace{5pt}
\scalebox{0.3}{\includegraphics*{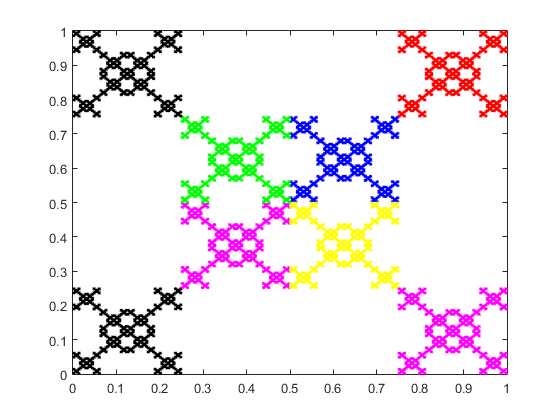}} \\
\caption{Fractals of Remark \ref{chain}}\label{fig2}
\end{figure}


\begin{remark}\label{chain} {\rm Note that if $D=\{d_1,...,d_q\}$ after re-indexing, it is possible that $d_i=d_j$ for $i\neq j$.
 Two examples are given in Figure \ref{fig2}.  For the first fractal of Figure \ref{fig2}, we have  $T=\tiny \begin{bmatrix}
  -3 & -1  \\
  0 & 3
\end{bmatrix},
$
 and  $D= \{ \tiny d_1= \left[ \begin{array}{cc} 0 \\
0 \end{array} \right],d_2=\left[ \begin{array}{cc} 2 \\
1 \end{array} \right], d_3=\left[ \begin{array}{cc} -1 \\
1 \end{array} \right], d_4=\left[ \begin{array}{cc} 1 \\
3 \end{array} \right], d_5=\left[ \begin{array}{cc} 2 \\
0 \end{array} \right], d_6=\left[ \begin{array}{cc} 2 \\
2 \end{array} \right] , d_7=\left[ \begin{array}{cc} -2 \\
2 \end{array} \right], d_8=\left[ \begin{array}{cc} -1 \\
3 \end{array} \right] \}.
$

 In view of Proposition \ref{neighbors1} and Proposition \ref{neighbors2},
Proposition \ref{connectedness} is quite feasible because the connectedness can be decided by a simple inspection of $D$ using
$\mathcal{N}$ in Propositions \ref{neighbors1}-\ref{neighbors2}. That is, we get a
graph-independent way of checking the connectedness.
For the second fractal on the right, $T= 4 I$ and
$D=\{ \tiny d_1=\left[
\begin{array}{cc}  0 \\ 0  \end{array}
\right], d_2=\left[
\begin{array}{cc}  1 \\ 1  \end{array}
\right], d_3=\left[
\begin{array}{cc}  2 \\ 2  \end{array}
\right],
 d_{4}=\left[
\begin{array}{cc}  3 \\ 3  \end{array}
\right],
 d_{5}=\left[
\begin{array}{cc}  2 \\ 1  \end{array}
\right],
 d_{6}=\left[
\begin{array}{cc}  1 \\ 2  \end{array}
\right],
 d_{7}=\left[
\begin{array}{cc}  0 \\ 3  \end{array}
\right],
 d_{8}=\left[
\begin{array}{cc}  3 \\ 0  \end{array}
\right]
 \}$. }

\end{remark}

For the general case $D\subset {\Bbb Z}^2$, we have the following trivial proposition, which again can be used together with Proposition \ref{connectedness}. Let
$$\mathcal{M}_1= \{ \pm (ke_1\pm le_2)  \ | \  \ k,l \in \mathbb{N} \ {\rm and} \ ka_1 \pm la_2\in  \Delta D\},$$
$$\mathcal{M}_0= \{ \pm (ke_1\pm le_2)  \ | \  \ k,l \in \mathbb{N} \ {\rm and} \ kb_1 \pm lb_2\in  \Delta D\}.$$
Note that it is possible that $\mathcal{M}_1=\emptyset$ or $\mathcal{M}_0=\emptyset$.

\begin{proposition}\label{generalneighbor1} \ Assume that $T$ is as in  (\ref{matrix}) with $t=1$, $n,m>0$  and $D\subset {\Bbb Z}^2$. Then
%
%
$\mathcal{N}\supseteq \{ \pm ke_i \ | \ k \in \mathbb{N}, \ i\in \{1,2 \} \ {\rm and} \ ka_i \in \Delta D \}\cup \mathcal{M}_1.$

%
%
%
%
\end{proposition}

\begin{proposition}\label{generalneighbor2} \ Assume that $T$ is as in  (\ref{matrix}) with $t=0$, $n,m>0$  and $D\subset {\Bbb Z}^2$. Let $b_2=a_2$, $b_3=a_3$.  Then

%
$\mathcal{N}\supseteq \{ \pm ke_i \ | \ k \in \mathbb{N}, \ i\in \{1,2,3 \} \ {\rm and} \ kb_i \in \Delta D \}\cup \mathcal{M}_0.$

%
%
%
\end{proposition}

\end{document}